\documentclass[12pt]{article}
\usepackage{amssymb, amsmath}

\newtheorem{theorem}{Theorem}[section]
\newtheorem{lemma}{Lemma}[section]
\newtheorem{definition}{Definition}[section]

\newcommand{\n}{\nonumber}

\newcommand{\si}{\sigma_R (|x|)}

\newcommand{\s}{\sigma}

\newcommand{\bb}{\begin{equation}}
\newcommand{\ee}{\end{equation}}
\newcommand{\bq}{\begin{eqnarray}}
\newcommand{\eq}{\end{eqnarray}}
\newcommand{\bqn}{\begin{eqnarray*}}
\newcommand{\eqn}{\end{eqnarray*}}

\begin{document}
\title{ On the nonexistence of  global weak solutions to the Navier-Stokes-Poisson equations in $\Bbb R^N$}
\author{ Dongho Chae\\
 Department of Mathematics\\
  Sungkyunkwan
University\\
 Suwon 440-746, Korea\\
 e-mail: {\it chae@skku.edu }}
 \date{}
\maketitle
\begin{abstract}
In this paper we prove nonexistence of stationary weak solutions to
the Euler-Poisson equations and the Navier-Stokes-Poisson equations
in $\Bbb R^N$, $N\geq 2$, under suitable  assumptions of
integrability for the density, velocity and the potential of the
force field. For the time dependent Euler-Poisson equations we prove
nonexistence result assuming additionally temporal asymptotic
behavior near infinity of the second moment of density. For a class
of time dependent Navier-Stokes-Poisson equations this asymptotic
behavior of the density can be proved if we assume the standard
energy inequality, and therefore the nonexistence of
global weak solution follows from more plausible assumption in this case.
{\bf keywords:} Euler-Poisson equations, Navier-Stokes Poisson
equations, nonexistence of global weak solutions
\end{abstract}
\section{Introduction}
 \setcounter{equation}{0}
We are concerned  on the  Navier-Stokes-Poisson equations in $\Bbb
R^N$, $N\geq 1$.
  $$
(NSP, EP)\left\{\aligned &\partial_t \rho + \mbox{div}(\rho v) = 0, \\
& \partial_t (\rho v) + \mbox{div}(\rho v \otimes v) = -\nabla
p+k\rho \nabla \Phi+ \mu \Delta v +(\mu +\lambda ) \nabla \mathrm{div }\, v, \\
&\Delta \Phi =\rho, \\
 &\rho\geq 0, p=p(\rho)\geq 0 (\mbox{$
p=0$ only if
 $\rho=0$}).\\
\endaligned \right.
$$
The system (NSP, EP) describes  compressible gas flows,  and $\rho,
v, \Phi$ and $p$ denote the density, velocity, the potential of the
underlying force and the pressure respectively. Here $k$ is a
physical constant, which signifies the property of the forcing,
repulsive if $k > 0$ and attractive if $k < 0$. {\it In this paper
we consider only the case of repulsive forcing, $k\geq 0$.}
 We treat the viscous
case $\mu
>0$, which corresponds to the Navier-Stokes-Poisson equations(NSP), and the inviscid case
$\mu=\lambda=0$, which corresponds to the Euler-Poisson
equations(EP), simultaneously.  Our aim here is to prove the
nonexistence of  global weak solutions to the system (NSP, EP) with
$N\geq 2$ under suitable integrability conditions for the solutions
together with additional condition for the initial data. For the
stationary case the standard finite energy condition already implies
the integrability for the solutions. For the time dependent
Euler-Poisson equations, however, we need an extra condition for the
temporal asymptotic behavior near infinity of the second moment of
density $\rho(x,t)$(see (\ref{5}) for more specification) to get
desired nonexistence result. For the Navier-Stokes-Poisson equations
describing the isothermal viscous fluids  with $p(\rho)=a
\rho^\gamma, 1<\gamma \leq N/4+1/2, N\geq 3$, the condition of the
asymptotic behavior of density can be proved, thanks to a lemma due
to Guo and Jiang(\cite{guo}), if we assume the energy inequality.
Hence, in this case the finite energy condition together with $v\in
L^{\frac{N}{N-1}} (\Bbb R^N\times[0,T))$ for all $T>0$ imply the
nonexistence of the global weak solutions satisfying the energy
inequality for  an initial data satisfying suitable sign condition.
This implies that even if the finite blow-up happens for certain
smooth initial data, it could not be continued as a physically
meaningful global weak solution afterwards. The results derived in
this paper for the nonexistence of global weak solutions could be
regarded as Luouville type of theorems. The convection term and the
forcing term have the ``positivity" structure, in appropriate sense,
which resembles the ellipticity in the elliptic partial differential
equations. Those ``positivity" structures, combined with the actual
positivity of the pressure term provides us the desired nonexistence
results for the nontrivial global weak solutions. Earlier
observations of similar feature for the convection term were made,
and applied to the
 compressible Euler and the compressible Navier-Stokes equations in
 \cite{cha1,cha2}, and the theorems obtained here are
 generalizations of those in \cite{cha1,cha2}, which are not
straightforward due to the nonlinear forcing term $k\rho\nabla\Phi$
in (NSP, EP). The sign condition for
 $k\geq 0$ and the restriction of spatial dimension $N\geq 2$ are
 crucially important to deduce the favorable positivity of the forcing
 term. At this moment we do not
 know if similar nonexistence results hold also for $k<0$ or $N=1$.

\subsection{Nonexistence of stationary weak solutions}
In this section we state  precisely the nonexistence theorem for the
stationary weak solutions to the system $(NSP,EP)$.
 A stationary weak
 solutions of  (NSP, EP)  is defined as follows.
  \begin{definition} We say that  a triple $(\rho, v,\Phi)\in  L^\infty_{loc}(\Bbb R^N )\times
  [L^2_{loc} (\Bbb R^N)]^N\times  W^{2,2}_{loc} (\Bbb R^N )$ is a  stationary weak solution of $(NSP, EP)$
  if
  \bq\label{11}
 && \int_{\Bbb R^N} \rho v\cdot \nabla \psi \,dx=0 \qquad \forall \psi
 \in C_0 ^\infty (\Bbb R^N),\\
 \label{12}
 &&\int_{\Bbb R^N} \rho v\otimes v :\nabla \phi \,dx=-\int_{\Bbb R^N} p\, \,\mathrm{div }\, \phi
 \,dx -k\int_{\Bbb R^N} \rho \nabla \Phi \cdot \phi \,dx-\mu\int_{\Bbb R^N} v\cdot \Delta \phi\, dx\n \\
 &&\hspace{1.in}\qquad -(\mu +\lambda
 )\int_{\Bbb R^N} v\cdot \nabla \mathrm{div}\, \phi \, dx
 \quad \forall \phi
 \in [C_0 ^\infty (\Bbb R^N)]^N,\n\\
 \ \\
 \label{12a}
 &&\Delta\Phi =\rho \quad\mbox{almost everywhere in $\Bbb R^N$},\\
 \label{14}
 &&p = p(\rho)\geq 0,\quad p=0 \,\,\mbox{only if}\,\, \rho=0\quad\mbox{almost everywhere on $\Bbb R^N$} .
  \eq
  \end{definition}
In our proof of the following theorem we do not use the equation of
continuity (\ref{11}).

 \begin{theorem} Let $N\geq 2$. Suppose  $(\rho, v, \Phi)$ is a
 stationary weak solution to $(NSP, EP)$ satisfying one of the following
 conditions depending on $\mu$ and $\lambda$.
 \begin{itemize}
 \item[(i)] For (EP)($\mu=\lambda =0$);
 \begin{description}
 \item[(i-a)] \rm{({\bf The case $N\geq 3$})} There exists $w\in L^1_{loc} ([0, \infty))$,
 which is positive, non-increasing on $[0, \infty)$ such that
\bq\label{15}
 \lefteqn{\int_{\Bbb R^N} (\rho |v|^2 +p+k|\nabla \Phi|^2)\times}\hspace{.0in}\n \\
 && \times\left[w(|x|) +\frac{1}{|x|}\int_0 ^{|x|} w(s)ds +\frac{1}{|x|^2}
 \int_0 ^{|x|} \int_0 ^r
 w(s)dsdr \right] dx
 <\infty.\n \\
 \eq
 \item[(i-b)] \rm{({\bf The case $N=2$})}
 \bb\label{15a}
\int_{\Bbb R^N} (\rho |v|^2 +p+k|\nabla \Phi|^2)dx <\infty.
 \ee
 \end{description}
 \item[(ii)] For (NSP) ($\mu >0$) with  $N\geq 2$ ;
 \begin{description}
 \item[(a)] if $2\mu+\lambda=0$,
 \bb\label{15b}
 \int_{\Bbb R^N} (\rho |v|^2 +p+k|\nabla \Phi|^2)\,dx
 <\infty.
\ee
 \item[(b)]
 if  $2\mu+\lambda\neq0$,
\bb\label{15c}
  \int_{\Bbb R^N} (\rho |v|^2
+|v|^{\frac{N}{N-1}}+p+k|\nabla \Phi|^2)\,dx
 <\infty.
 \ee
 \end{description}
 \end{itemize}
 Then, $\rho(x)=0, \nabla\Phi (x)=0$ for almost every $x\in \Bbb R^N$.
\end{theorem}
{\it Remark 1.1 } Choosing, in particular,
 \bb\label{wa}
w(r)=1/(1+r^2),
 \ee
 then for all $x\in \Bbb R^N$ we have
 \bqn
&&\int_0 ^{|x|} w(r)dr \leq \frac{\pi}{2},\qquad\int_0 ^{|x|} \int_0 ^r w(s)ds dr \leq \frac{\pi|x|}{2},\quad \mbox{and}\\
&& w(|x|) +\frac{1}{|x|}\int_0 ^{|x|} w(s)ds +\frac{1}{|x|^2}
 \int_0 ^{|x|} \int_0 ^r
 w(s)dsdr \leq \frac{C}{1+|x|}
 \eqn
 for some constant $C$ independent of $x$.
Thus the condition for the initial data (\ref{4a}) and (\ref{4b})
are implied by
 \bb\label{4ca}\int_{\Bbb R^N} \rho_0 (x)|v_0 (x)|dx <\infty,
 \ee
 and
\bb\label{4caa}
 \int_{\Bbb R^N} \rho_0 (x)v_0 (x)\cdot \frac{x}{|x|}
\arctan(|x|) dx \geq CK_1
 \ee
 respectively, while the condition (\ref{15}) is
 implied by
 \bb\label{4ea}
 \int_{\Bbb R^N} \frac{\rho(x)|v(x)|^2 +p(x)+k|\nabla \Phi (x)|^2}{1+|x|}dx
 <\infty
 \ee
 respectively.
Note that the condition (\ref{4ea}) is even  weaker than the finite
energy condition, in the sense that it is implied by the finite
energy condition that is also obtained from (\ref{15}) by choosing
$w=1$.

\subsection{Nonexistence of time dependent weak solutions}
  The definition of time dependent weak solutions for (NSP, EP) is follows.
  \begin{definition} We say a triple
  $$(\rho,v, \Phi)\in  L^1_{loc} ((0, \infty);L^\infty_{loc}(\Bbb R^N ))\times
  [L^1_{loc} ((0, \infty); L^2_{loc} (\Bbb R^N))]^N\times
L^2_{loc}((0,\infty); W^{2,2}_{loc} (\Bbb R^N))
  $$
   is a global weak solution of
  (NSP) with initial data $(\rho_0, v_0)$
  if
  \bq\label{1}
 && \xi(0)\int_{\Bbb R^N} \rho_0 (x) \psi (x)dx+
 \int_0 ^\infty \int_{\Bbb R^N} \rho (x,t) \psi (x) \xi ' (t)dxdt\n\\
 &&\qquad+\int_0 ^\infty\int_{\Bbb R^N} \rho v(x,t)\cdot \nabla \psi (x)\xi(t)\,dx=0
 \n \\
 &&\hspace{2.in}\forall \psi
 \in C_0 ^\infty (\Bbb R^N),\xi \in C_0 ^1([0, \infty)),
\eq
 \bqn
 &&\xi(0)\int_{\Bbb R^N} \rho_0 (x)v_0(x)\cdot \phi (x)dx+
 \int_0 ^\infty \int_{\Bbb R^N} \rho (x,t)v(x,t)\cdot \phi (x) \xi ' (t)dxdt\n \\
 &&\qquad+\int_0 ^\infty\int_{\Bbb R^N} \rho(x,t) v(x,t)\otimes v(x,t) :\nabla \phi (x)\xi(t)\,dxdt
 \eqn
 \bq
 \label{2}
 &&\qquad=
 -\int_0 ^\infty\int_{\Bbb R^N} p(x,t)\, \,\mathrm{div }\,
 \phi(x)\xi(t)
 \,dxdt -k\int_0 ^\infty \int_{\Bbb R^N} \rho (x,t) \nabla\Phi(x,t)\cdot \phi(x)
 \xi(t)dxdt\n \\
 &&\,\qquad-\mu\int_0 ^\infty\int_{\Bbb R^N} v(x,t)\cdot \Delta \phi(x)\xi(t)\, dxdt-(\mu +\lambda
 )\int_0 ^\infty\int_{\Bbb R^N} v(x,t)\cdot \nabla \mathrm{div}\, \phi (x)\xi(t)\,
 dxdt\n \\
 && \hspace{2.in}\qquad\forall \phi
 \in [C_0 ^\infty (\Bbb R^N)]^N,\xi \in C_0 ^1([0, \infty)),\\
 \label{2a}
  &&\Delta\Phi =\rho \quad\mbox{almost everywhere on $\Bbb R^N\times[0, \infty )$},\\
 \label{4}
 &&\rho\geq 0, p=p(\rho)\geq 0 (\mbox{$ p=0$ only if
 $\rho=0$})\quad\mbox{almost everywhere on $\Bbb R^N\times[0, \infty )$}.
  \eq
  \end{definition}
 In the above the derivatives of $\xi \in C^1_0([0, \infty))$ at
 $t=0$ should be understood as $\xi'(0):=\xi'(0+)$.
 \begin{theorem}[Conditional nonexistence for (EP)]
 \begin{itemize}
 \item[(i)] \underline{The case $N\geq 3$ }:
 Let the function  $w\in L^1_{loc} ([0, \infty))$ be given,
 which is positive, non-increasing on $[0, \infty)$, and
let $(\rho_0, v_0)$ satisfy
 \bb\label{4a}
 \int_{\Bbb R^N} \rho_0 (x)|v_0 (x)| \left[\int_0 ^{|x|} w(r)dr\right]dx
 <\infty.
 \ee
Suppose  $(\rho, v, \Phi)$ is a global weak solution to (EP) with
the initial data $(\rho_0, v_0)$ such that
 \bb\label{5}
 \lim\sup_{\tau \to \infty}
 \int_{\tau \leq t\leq 2\tau}\int_{\Bbb R^N} \frac{\rho (x,t)}{1+t^2} \left[\int_0 ^{|x|} \int_0 ^r w(s)ds
 dr\right]
 dxdt \leq K_1
 \ee
 for a constant $K_1\geq 0$,
satisfying
  \bq\label{6}
 \lefteqn{\int_0 ^T\int_{\Bbb R^N} (\rho |v|^2 +p+k|\nabla \Phi |^2)\times}\hspace{.0in}\n \\
 && \times\left[w(|x|) +\frac{1}{|x|}\int_0 ^{|x|} w(s)ds +\frac{1}{|x|^2}
 \int_0 ^{|x|} \int_0 ^r
 w(s)dsdr \right] dxdt
 <\infty\n \\
 \eq
 for all $T>0$.
Then, necessarily the following inequality holds true.
  \bqn
 \lefteqn{\int_{0} ^\infty\int_{\Bbb R^N} \rho(x,t)\left[
w(|x|) \frac{(v\cdot x)^2}{|x|^2} +\frac{1}{|x|} \int_0 ^{|x|}
w(r)dr \left(|v|^2-\frac{(v\cdot x)^2}{|x|^2}\right) \right]
 \,dxdt}\n \\
&&\qquad+\int_{0} ^\infty\int_{\Bbb R^N}p(x,t)\left[w
 (|x|) +\frac{N-1}{|x|}\int_0 ^{|x|} w(r)dr\right] \, dxdt
    \eqn
\bq\label{6a}
 &&\,\qquad+\frac{\left(
 N-3\right)k}{2}\int_{0} ^\infty\int_{\Bbb R^N}\frac{|\nabla \Phi
 |^2}{|x|}\left[\int_0 ^{|x|}w (s) ds\right]\,dxdt\n \\
 &&\,\qquad+k\int_{0} ^\infty\int_{\Bbb R^N}
 \left[\frac{1}{|x|}\int_0 ^{|x|}w (s) ds-w(|x|)
\right]\frac{(x\cdot \nabla \Phi )^2}{|x|^2}\,dxdt\n \\
 &&\qquad+\frac{k}{2}\int_{0} ^\infty\int_{\Bbb R^N} |\nabla \Phi |^2  w
 (|x|)\,dxdt\n \\
 &&\qquad+\int_{\Bbb R^N} \rho_0 (x)v_0(x)\cdot
\frac{x}{|x|}\left[ \int_0 ^{|x|} w(r)dr\right]\,dx\leq C K_1
 \eq
 for a constant $C$.
 Therefore, if
  \bb\label{4b}
 \int_{\Bbb R^N} \rho_0 (x)v_0 (x)\cdot \frac{x}{|x|} \left[\int_0 ^{|x|} w(r)dr\right] dx
 >
CK_1,
 \ee
then there exists no  global weak solution satisfying
(\ref{5})-(\ref{6}).
\item[(ii)] \underline{The case $N=2$ :}
Let $(\rho_0, v_0)$ satisfy
 \bb\label{4aa}
 \int_{\Bbb R^N} \rho_0 (x)|v_0 (x)||x|dx
 <\infty.
 \ee
Suppose  $(\rho, v, \Phi)$ is a global weak solution to (EP) with
the initial data $(\rho_0, v_0)$ such that
 \bb\label{5a}
 \lim\sup_{\tau \to \infty}
 \int_{\tau \leq t\leq 2\tau}\int_{\Bbb R^N} \frac{\rho
 (x,t)}{1+t^2}|x|^2
 dxdt \leq K_1
 \ee
 for a constant $K_1\geq 0$,
satisfying
  \bb\label{6a}
\int_0 ^T\int_{\Bbb R^N} (\rho |v|^2 +p+k|\nabla \Phi |^2)dxdt
 <\infty
 \ee
 for all $T>0$.
Then, necessarily the following inequality holds true.
  \bb\label{6aa}
 \int_{0} ^\infty\int_{\Bbb R^2} (\rho|v|^2+
2p) dxdt +\int_{\Bbb R^2} \rho_0 (x)v_0(x)\cdot x\,dx\leq C K_1
 \ee
 for a constant $C$.
 Therefore, if
  \bb\label{4ba}
 \int_{\Bbb R^2} \rho_0 (x)v_0 (x)\cdot x \,dx
 >
CK_1,
 \ee
then there exists no  global weak solution satisfying
(\ref{5a})-(\ref{6aa}).
\end{itemize}
\end{theorem}
{\it Remark 1.2 } Similarly to Remark 1.1, choosing $
w(r)=1/(1+r^2)$, the conditions for the solution (\ref{5}) and
(\ref{6}) are
 implied by
 \bb\label{4d}
 \lim\sup_{\tau \to \infty}\int_{\tau\leq t\leq2\tau} \int_{\Bbb R^N}\frac{\rho (x,t)|x|}{1+t^2}
 dxdt \leq K_1,
  \ee
and
 \bb\label{4e}
 \int_0 ^T\int_{\Bbb R^N} \frac{\rho(x,t)|v(x,t)|^2 +p(x,t)+k|\nabla \Phi (x,t)|^2}{1+|x|}dxdt
 <\infty\qquad\forall T>0
 \ee
 respectively.\\
 \ \\
As will be seen below, for a class of isothermal viscous fluids the
key condition (\ref{5}) is really satisfied with $K_1=0$, if we
assume the energy inequality, and we have the following stronger
nonexistence results of the global weak solutions.

\begin{theorem}[Nonexistence for (NSP)]
We fix $N\geq 3$, $1<\gamma \leq N/4 +1/2$,  $\mu>0 , \mu+\lambda
>0$, and the following form of
pressure law,
 \bb\label{pre}
 p=p(\rho)=a \rho^\gamma
 \ee
 in (NSP).
 Let the initial
 data $(\rho_0, v_0)$ satisfy
 \bb\label{7}
 \int_{\Bbb R^N} \rho_0 (x)|v_0 (x)| |x|dx <\infty.
 \ee
Suppose $(\rho, v, \Phi)$  is a global weak solution to (NSP) such
that
  \bb\label{11aa}
  \int_0 ^T\int_{\Bbb R^N} \left[\rho |v|^2+p
+|v|^{\frac{N}{N-1}} +k|\nabla \Phi |^2\right]\,dxdt
 <\infty
 \ee
 for all $T>0$.
We further assume  that the following energy inequality holds.
 \bq\label{energy}
&&E(t)+\int_0 ^t \int_{\Bbb R^N} (\mu |\nabla v|^2 +(\mu+\lambda) |
\mathrm{div}\, v|^2 )
dxds\leq E(0)<\infty \quad \forall t\geq 0, \n \\
&&\quad\mbox{where} \qquad E(t):=\int_{\Bbb R^N} \left[\frac
{\rho}{2} |v|^2 +\frac{a\rho^{\gamma} }{\gamma-1}+\frac{k}{2}
|\nabla \Phi |^2\right] dx. \eq Then, necessarily we have the
equality
  \bq\label{11a}
 &&\int_{0} ^\infty\int_{\Bbb R^N} \left[\rho(x,t) |v(x,t)|^2 +N p(x,t) +\frac{N-2}{2}
 |\nabla \Phi (x,t)|^2\right]\, dxdt\n \\
&&\qquad=-\int_{\Bbb R^N} \rho_0 (x)v_0(x)\cdot x\,dx
 \eq
 Therefore, if
  \bb\label{11b}
 \int_{\Bbb R^N} \rho_0 (x)v_0 (x)\cdot x\, dx \geq
 0,
 \ee
then the only global weak solution
 corresponds to $\rho=0$ and $\nabla \Phi=0$ almost everywhere on $\Bbb R^N\times [0, \infty)$. In particular, if the strict
 inequality holds in (\ref{11b}), then there exists no global weak solution
 satisfying (\ref{7})-(\ref{energy}).
\end{theorem}

\section{Proof of the main theorems}
\setcounter{equation}{0}

\noindent{\bf Proof of Theorem 1.1 } We suppose there exists a
stationary weak solution $(\rho, v, \Phi)$. We begin the proof with
the inviscid case.\\
\ \\
 \noindent{\em \underline{(i) The case
$\mu=\lambda=0$ (EP): }}
 Let us
consider a radial cut-off function $\sigma\in C_0 ^\infty(\Bbb R^N)$
such that
 \bb\label{18}
   \sigma(|x|)=\left\{ \aligned
                  &1 \quad\mbox{if $|x|<1$}\\
                     &0 \quad\mbox{if $|x|>2$},
                      \endaligned \right.
 \ee
and $0\leq \sigma  (x)\leq 1$ for $1<|x|<2$. We set
 \bb\label{19}
W(u):=\int_0 ^{u} \int_0 ^s w(r)drds.
 \ee
Then, for each $R >0$, we define
 \bb\label{110}
\varphi_R (x)=W(|x|)\s \left(\frac{|x|}{R}\right)=W(|x|)\s_R
(|x|)\in C_0 ^\infty (\Bbb R^N).
 \ee
We   choose the vector test function
 $\phi$ in (\ref{12}) as
 \bb\label{111}
  \phi= \nabla \varphi_R (x).
 \ee
 Then, after routine computations, the equation (\ref{12}) becomes
 \bq\label{112}
 \lefteqn{0=\int_{\Bbb R^N} \rho (x) \left[W^{\prime\prime} (|x|)
\frac{(v\cdot x)^2}{|x|^2} + W^{\prime} (|x|)
\left(\frac{|v|^2}{|x|}
-\frac{(v\cdot x)^2}{|x|^3}\right) \right] \si  \,dx }\hspace{.5in}\n\\
&&\quad+\int_{\Bbb R^N} \rho(x) W' (|x|) \s'
\left(\frac{|x|}{R}\right)
\frac{(v\cdot x)^2}{R|x|^2} \,dx\n \\
&&\quad+ \frac{1}{R}\int_{\Bbb R^N} \rho(x)\left( \frac{|v|^2}{|x|}
-\frac{(v\cdot x)^2}{|x|^3} \right)
\s'\left(\frac{|x|}{R}\right)W(|x|) \,dx\n \\
  &&\quad+\int_{\Bbb R^N} \rho(x)\frac{(v\cdot x)^2}{
R^2|x|^2} \s^{\prime\prime} \left(\frac{|x|}{R}\right) W(|x|)
\,dx \n \\
 &&\quad+ \int_{\Bbb R^N}p(x)\left[ W^{\prime\prime}
 (|x|) +(N-1)\frac{W' (|x|)}{|x|}\right]\sigma_R (|x|)
  \, dx\n \\
&&\quad+ \frac{2}{R}\int_{\Bbb R^N}p(x) W' (|x|)
 \s' \left(\frac{|x|}{R}\right)\, dx\n \\
&&\quad+ \frac{N-1}{R}\int_{\Bbb R^N}p(x)\frac{1}{|x|}\s'
\left(\frac{|x|}{R}\right) W(|x|) \, dx\n \\
&&\quad+ \int_{\Bbb R^N}p(x)\frac{1}{R^2} \s^{\prime\prime}
\left(\frac{|x|}{R}\right)W(|x|) \,
dx \n \\
&&\quad +k\int_{\Bbb R^N} \rho \nabla \Phi \cdot \nabla[ W(|x|)\si ]\,dx \n \\
 &&:=I_1+\cdots +I_9.
\eq
 In terms of the function $W$ defined in  (\ref{19}) our
condition (\ref{15})
 can be written as
\bb\label{114} \int_{\Bbb R^N} (\rho(x)|v(x)|^2 +|p(x)| +k|\nabla
\Phi |^2) \left[W^{\prime\prime}
 (|x|) +\frac{1}{|x|}W' (|x|) +\frac{1}{|x|^2} W(|x|)\right]
 dx<\infty.
 \ee
Since \bqn
 &&\int_{\Bbb R^N} \rho(x)\left|\left[W^{\prime\prime} (|x|)
\frac{(v\cdot x)^2}{|x|^2} + W^{\prime} (|x|)
\left(\frac{|v|^2}{|x|} -\frac{(v\cdot x)^2}{|x|^3}\right)
\right]\right|dx\\
&&\qquad \leq 2
  \int_{\Bbb R^N} \rho(x)|v(x)|^2\left[ W^{\prime\prime}
(|x|)+\frac{W'(|x|)}{|x|} \right] dx <\infty, \eqn
  we can use the dominated convergence theorem to show that
  \bb\label{115}
  I_1 \to \int_{\Bbb R^N} \rho(x)\left[W^{\prime\prime} (|x|)
\frac{(v\cdot x)^2}{|x|^2} + W^{\prime} (|x|)
\left(\frac{|v|^2}{|x|} -\frac{(v\cdot x)^2}{|x|^3}\right) \right]
\,dx
 \ee
  as $R\to \infty$.
Similarly,
  \bb\label{116}
  I_5\to \int_{\Bbb R^N}p(x)\left[ W^{\prime\prime}
 (|x|) +(N-1)\frac{W' (|x|)}{|x|}\right]
 \, dx
\ee as $R\to \infty$.
 For $I_2$ we estimate
 \bq\label{117}
  |I_2 |&\leq& \int_{R<|x|<2R} \rho(x)|v(x)|^2\left|\s'
\left(\frac{|x|}{R}\right)\right|
  \frac{W'(|x|)}{|x|} \frac{|x|}{R}dx\n \\
  &\leq &2 \sup_{1<s<2} |\s'(s)|
\int_{R<|x|<2R}\rho(x) |v(x)|^2
  \frac{W'(|x|)}{|x|} dx\n \\
\to 0
  \eq
 as $R\to \infty$ by the dominated convergence theorem.
Similarly
  \bq\label{118}
   |I_3|&\leq &2 \int_{R<|x|<2R}  \frac{|x|}{R} \rho(x) |v(x)|^2
\left|\s'\left(\frac{|x|}{R}\right)\right|\frac{W(|x|)}{|x|^2} \,dx \n \\
&\leq &4\sup_{1<s<2} |\s'(s)|
 \int_{R<|x|<2R}\rho(x) |v(x)|^2
  \frac{W'(|x|)}{|x|}dx
\to 0,\n \\
  \eq
  and
  \bq\label{119}
  |I_4|&\leq&\int_{R<|x|<2R}\frac{|x|^2}{R^2} \rho(x)|v(x)|^2\left|\s^{\prime\prime}
  \left(\frac{|x|}{R}\right)\right|\frac{W(|x|)}{|x|^2} \,dx\n \\
  &\leq&4\sup_{1<s<2} |\s^{\prime\prime}(s)|
 \int_{R<|x|<2R}\rho(x)|v(x)|^2
  \frac{W(|x|)}{|x|^2}\,dx \to 0\n \\
  \eq
   as $R\to \infty$. The estimates for $I_6,I_7$ and $I_8$ are
   similar to the above, and we find
   \bq\label{120}
   |I_6|&\leq &2 \int_{R<|x|<2R}|p(x)| \frac{|x|}{R}\frac{W'
   (|x|)}{|x|}
 \left|\s' \left(\frac{|x|}{R}\right)\right|  \, dx\n \\
 &\leq& 4\sup_{1<s<2} |\s' (s)|
\int_{R<|x|<2R}|p(x)|\frac{W'
   (|x|)}{|x|}dx \to 0,\n \\
  \eq
\bq\label{121}
   |I_7|&\leq & (N-1)\int_{R<|x|<2R}|p(x)|\frac{|x|}{R}\left|\s'
\left(\frac{|x|}{R}\right)\right| \frac{W(|x|)}{|x|^2} \, dx\n \\
 &\leq& 2\sup_{1<s<2} |\s' (s)|
\int_{R<|x|<2R}|p(x)|\frac{W(|x|)}{|x|^2}dx\to 0,
  \eq
  and
  \bq\label{122}
 |I_8|&\leq&\int_{\Bbb R^N}|p(x)|\frac{|x|^2}{R^2}
\left|\s^{\prime\prime}\left(\frac{|x|}{R}\right)\right|
\frac{W(|x|)}{|x|^2} \, dx \n\\
  &\leq& 4\sup_{1<s<2}
|\s^{\prime\prime} (s)|
 \int_{R<|x|<2R}|p(x)|\frac{W
   (|x|)}{|x|^2}dx \to 0
  \eq
  as $R\to \infty$ respectively. Using the relation (\ref{2a}), and integrating by parts we
  compute
 \bq\label{122a}
 I_9&=&k\int_{\Bbb R^N} \Delta \Phi \nabla \Phi \cdot \nabla
 [ W(|x|) \si ]\, dx\n \\
 &=&-k\sum_{i,j=1}^N \int_{\Bbb R^N} \partial_i \Phi \partial_j\Phi
 \partial_i\partial_j
 [ W(|x|) \si ]\, dx -k\sum_{i,j=1}^N \int_{\Bbb R^N} \partial_i \Phi \partial_j\partial_i\Phi
\partial_j
 [ W(|x|) \si ]\, dx\n \\
 &=&-k\sum_{i,j=1}^N \int_{\Bbb R^N} \partial_i \Phi \partial_j\Phi
 \partial_i\partial_j
 [ W(|x|) \si ]\, dx +\frac{k}{2}  \int_{\Bbb R^N} |\nabla
 \Phi|^2\Delta
 [ W(|x|) \si ]\, dx\n \\
&=&-k\int_{\Bbb R^N} \left[\frac{|\nabla \Phi |^2}{|x|}
-\frac{(x\cdot \nabla \Phi )^2}{|x|^3} \right]\left[
W'(|x|)\si + \frac{W(|x|)}{R}\sigma' \left(\frac{|x|}{R}\right)\right]\,dx\n \\
 &&\,-k\int_{\Bbb R^N} \frac{(x\cdot \nabla \Phi
)^2}{|x|^2} \left[ W^{\prime\prime} (|x|)\si
+\frac{W(|x|)}{R^2}\sigma^{\prime\prime}\left(\frac{|x|}{R}\right)
+\frac{2 W^{\prime}
(|x|)}{R}\sigma^{\prime}\left(\frac{|x|}{R}\right)\right]\,dx\n \\
&&+\frac{(N-1)k}{2}\int_{\Bbb R^N} \frac{|\nabla \Phi |^2}{|x|}
\left[
W'(|x|)\si + \frac{W(|x|)}{R}\sigma' \left(\frac{|x|}{R}\right)\right]\,dx\n \\
  &&+\frac{k}{2}\int_{\Bbb R^N} |\nabla \Phi |^2 \left[ W^{\prime\prime} (|x|)\si
+\frac{W(|x|)}{R^2}\sigma^{\prime\prime}\left(\frac{|x|}{R}\right)
+\frac{2 W^{\prime}
(|x|)}{R}\sigma^{\prime}\left(\frac{|x|}{R}\right)\right]\,dx.\n \\
&:=&J_1+\cdots+J_4.
 \eq
 By
similar computations to (\ref{115})-(\ref{122}), using (\ref{114}),
we find that
 \bqn
   J_1&=&-k\int_{\Bbb R^N} \left[\frac{|\nabla \Phi
|^2}{|x|} -\frac{(x\cdot \nabla \Phi )^2}{|x|^3} \right]W'(|x|)\,dx
+o(1),\\
J_2&=&-k\int_{\Bbb R^N} \frac{(x\cdot \nabla \Phi )^2}{|x|^2}
W^{\prime\prime} (|x|)\,dx+o(1), \eqn \bqn
J_3&=&\frac{(N-1)k}{2}\int_{\Bbb R^N} \frac{|\nabla \Phi |^2}{|x|}
 W'(|x|)\, dx +o(1),\n \\
 J_4&=&\frac{k}{2}\int_{\Bbb R^N} |\nabla \Phi |^2  W^{\prime\prime}
 (|x|)\,dx+o(1)
 \eqn
as $R\to \infty$. Therefore, taking the limit $R\to \infty$, and
rearranging the remaining terms, we have
 \bq\label{122aa}
 I_9&\to &-k\int_{\Bbb R^N} \left[\frac{|\nabla \Phi |^2}{|x|}
-\frac{(x\cdot \nabla \Phi )^2}{|x|^3} \right]W'(|x|)\,dx
-k\int_{\Bbb R^N} \frac{(x\cdot \nabla \Phi )^2}{|x|^2}
W^{\prime\prime} (|x|)\,dx\n \\
&&+\frac{(N-1)k}{2}\int_{\Bbb R^N} \frac{|\nabla \Phi |^2}{|x|}
 W'(|x|)\, dx
 +\frac{k}{2}\int_{\Bbb R^N} |\nabla \Phi |^2  W^{\prime\prime}
 (|x|)\,dx\n \\
 &&=\frac{\left(
 N-3\right)k}{2}\int_{\Bbb R^N}\frac{|\nabla \Phi
 |^2}{|x|}W'(|x|)\,dx+k\int_{\Bbb R^N}\left[\frac{W'(|x|)}{|x|}-W^{\prime\prime}
 (|x|)\right]\frac{(x\cdot \nabla \Phi )^2}{|x|^2}\,dx\n \\
 &&\qquad+\frac{k}{2}\int_{\Bbb R^N} |\nabla \Phi |^2  W^{\prime\prime}
 (|x|)\,dx
 \eq
  as $R\to \infty$. Thus passing $R\to \infty$ in (\ref{112}), and
  using (\ref{115})-(\ref{122aa}),
  we finally obtain
\bq\label{123}
 &&\int_{\Bbb R^N} \rho(x)\left[W^{\prime\prime} (|x|)
\frac{(v\cdot x)^2}{|x|^2} + W^{\prime} (|x|)
\left(\frac{|v|^2}{|x|} -\frac{(v\cdot x)^2}{|x|^3}\right) \right]
 \,dx \n \\
&&\qquad+\int_{\Bbb R^N}p(x)\left[ W^{\prime\prime}
 (|x|) +(N-1)\frac{W' (|x|)}{|x|}\right] \, dx\n \\
 &&\qquad+\frac{\left(
 N-3\right)k}{2}\int_{\Bbb R^N}\frac{|\nabla \Phi
 |^2}{|x|}W'(|x|)\,dx+k\int_{\Bbb R^N}\left[\frac{W'(|x|)}{|x|}-W^{\prime\prime}
 (|x|)\right]\frac{(x\cdot \nabla \Phi )^2}{|x|^2}\,dx\n \\
 &&\qquad+\frac{k}{2}\int_{\Bbb R^N} |\nabla \Phi |^2  W^{\prime\prime}
 (|x|)\,dx=0,
\eq
 which can be written, in terms of
the function $w(r)$,  as
 \bqn
 &&\int_{\Bbb R^N} \rho(x) \left[w (|x|)
\frac{(v\cdot x)^2}{|x|^2} + \frac{1}{|x|}\int_0 ^{|x|}w (s) ds
\left(|v|^2 -\frac{(v\cdot x)^2}{|x|^2}\right)\right]
\,dx \n \\
&&\qquad+\int_{\Bbb R^N}p(x)\left[ w
 (|x|) +\frac{N-1}{|x|}\int_0 ^{|x|} w (s)ds \right]\,dx\n \\
 &&\,\qquad+\frac{\left(
 N-3\right)k}{2}\int_{\Bbb R^N}\frac{|\nabla \Phi
 |^2}{|x|}\left[\int_0 ^{|x|}w (s) ds\right]\,dx
 \eqn
 \bq\label{124}
 &&\,\qquad+k\int_{\Bbb R^N}
 \left[\frac{1}{|x|}\int_0 ^{|x|}w (s) ds-w(|x|)
\right]\frac{(x\cdot \nabla \Phi )^2}{|x|^2}\,dx\n \\
 &&\quad\qquad+\frac{k}{2}\int_{\Bbb R^N} |\nabla \Phi |^2  w
 (|x|)\,dx=0.
 \eq
We note that
$$w (|x|)
\frac{(v\cdot x)^2}{|x|^2} + \frac{1}{|x|}\int_0 ^{|x|}w (s) ds
\left(|v|^2 -\frac{(v\cdot x)^2}{|x|^2}\right)\geq 0, $$
 and
 $$
 w
 (|x|) +\frac{N-1}{|x|}\int_0 ^{|x|} w (s)ds>0.
 $$
 Moreover, since $w(r)$ is non-increasing a.e. on $[0, \infty)$ by hypothesis, we have
 $$
 \frac{1}{|x|}\int_0 ^{|x|}w (s) ds-w(|x|)\geq 0\quad \mbox{for almost
 every $x\in \Bbb R^N$}.
 $$
 Thus all of the terms  in (\ref{124}) are nonnegative for $N\geq 3$, and we need to have
 $$
 p(x)=p(\rho(x))=0, \nabla \Phi(x)=0\qquad\mbox{almost every $x\in \Bbb
 R^N$}.
$$
Therefore $\rho(x)=0, \nabla\Phi(x)=0$ for almost every $x\in \Bbb R^N$.\\
 If $N=2$, then we fix $w(r)=1$ on $[0, \infty)$ in all of the computations
 leading to (\ref{124}). Then (\ref{124}) reduces to
\bb\label{124a}
 \int_{\Bbb R^2} \left[\rho(x)|v(x)|^2 +2p(x)\right] \,dx=0,
 \ee
 from which we have $\rho=0$  on $\Bbb R^2$. From
 (\ref{12a}) $\nabla \Phi$ is harmonic in $\Bbb R^N$, and this combined
 with the condition (\ref{15a}) implies $\nabla \Phi=0$ on $\Bbb
 R^N$.\\
 \ \\
 \noindent{\underline{\em (ii) The case of $\mu >0, N\geq 2$ with either $2\mu +\lambda =0$ or  $2\mu +\lambda \neq0$ (NSP): }}
 In this case
 we choose the  function $w(r)\equiv
 1$ on $[0, \infty)$ in the proof of (i) above, which is equivalent
 to the choice of the vector test function,
 $\phi=\frac12\nabla [|x|^2\si ]$ instead of (\ref{111}).
  We just need to show
 the vanishing of the viscosity term
 \bb\label{125}
 \mu\int_{\Bbb R^N} v\cdot \Delta \phi dx +(\mu +\lambda)\int_{\Bbb
 R^N} v\cdot \nabla \mathrm{div}\, v dx=o(1)
 \ee
 as $R\to \infty$.\\
 If $2\mu+\lambda =0$, then
\bqn
  J&:=&\mu\int_{\Bbb R^N} v\cdot \Delta \nabla (|x|^2 \s_R )\,
dx +(\mu +\lambda)\int_{\Bbb
 R^N} v\cdot \nabla [\mathrm{div}\, \nabla (|x|^2 \s_R)]\, dx\\
 &=&(2\mu +\lambda)\int_{\Bbb R^N} v\cdot\nabla \Delta (|x|^2
 \s \left(\frac{|x|}{R}\right)\, dx=0,
 \eqn
and (\ref{125}) holds true. \\
If $2\mu+\lambda \neq0$, then we compute and estimate
 \bqn |J|&=&|2\mu +\lambda| \left|\int_{\Bbb R^N} v\cdot\nabla \Delta (|x|^2
 \s \left(\frac{|x|}{R}\right)\, dx\right|\\
&\leq& |2\mu +\lambda| \left|\int_{\Bbb R^N}
(N+5)\left[\frac{(v\cdot x)}{R|x|}\s ' \left(\frac{|x|}{R}\right)+
\frac{(v\cdot x)}{R^2}\s^{\prime\prime}
\left(\frac{|x|}{R}\right)\right]\, dx\right|\n \\
&&\qquad+|2\mu +\lambda|\left|\int_{\Bbb R^N} \frac{|x|(v\cdot
x)}{R^3}\s^{\prime\prime\prime}
\left(\frac{|x|}{R}\right)\,dx \right|\n \\
&\leq&  \frac{C}{R} \int_{R\leq |x|\leq 2R} |v(x)| \, dx \leq C
\left(\int_{R\leq  |x|\leq 2R} |v(x)|^{\frac{N}{N-1}}
dx\right)^{\frac{N-1}{N}} \to 0
 \eqn
   as $R\to \infty$, since $v\in L^{\frac{N}{N-1}} (\Bbb R^N)$ by the hypothesis in the viscous case.
    Thus (\ref{125})  holds true.
 $\square$\\
 \ \\
\noindent{\bf Proof of Theorem 1.2 } Suppose there exists a global
weak solution $(\rho, v, \Phi)$ satisfying (\ref{1})-(\ref{4}) with
$\mu=\lambda=0$. We choose vector test function
$$\phi =\nabla \varphi_R (x),
$$
where $\varphi_R$ is defined in (\ref{18})-(\ref{110}).
 We also introduce $\eta \in C^\infty_0 ([0, \infty))$ as follows.
\bb\label{z15}
   \eta(t)=\left\{ \aligned
                  &1 \quad\mbox{if $0\leq t<1$}\\
                     &0 \quad\mbox{if $t>2$},
                      \endaligned \right.
 \ee
 and $0\leq \eta (t)\leq 1$ for all $t\geq 0$.
Then, we set
 \bb\label{z16}
 \eta_\tau(t)=\eta\left(\frac{t}{\tau}\right).
 \ee
We substituting $\phi(x)=\nabla \varphi_R (x), \xi(t)=\eta_\tau(t)$
into (\ref{2}). Then, substituting $\rho=\Delta \Phi$, and following
similar the computations  in (\ref{122a}), we obtain
  \bqn
  \lefteqn{0=\int_{\Bbb R^N} \rho_0 (x)v_0(x)\cdot \frac{x}{|x|} W'(|x|) \si dx}\n \\
  &&+\frac{1}{ R}\int_{\Bbb R^N}
  \rho_0 (x)v_0(x)\cdot \frac{x}{ |x|}W(|x|)\s'\left(\frac{|x|}{R}\right)
 dx\n \\
 &&+\int_{0} ^\infty\int_{\Bbb R^N} \rho (x,t)v(x,t) \cdot
 \nabla \varphi_R (x) \eta_\tau ^{\prime}(t) dxdt\n \\
&&+\int_{0} ^\infty\int_{\Bbb R^N} \rho (x,t) \left[W^{\prime\prime}
(|x|) \frac{(v\cdot x)^2}{|x|^2} +\right. \n\\
 &&\hspace{.5in}\left.+ W^{\prime} (|x|) \left(\frac{|v(x,t)|^2}{|x|}
-\frac{(v(x,t)\cdot x)^2}{|x|^3}\right) \right] \si \eta_\tau(t)
\,dx dt \n \\
  &&+\frac{1}{R}\int_{0} ^\infty\int_{\Bbb
R^N} \rho(x,t) W' (|x|) \s' \left(\frac{|x|}{R}\right)
\frac{(v(x,t)\cdot x)^2}{|x|^2}\eta_\tau(t) \,dxdt \n \\
 &&+ \frac{1}{R}\int_{0} ^\infty\int_{\Bbb R^N}
\rho(x,t)\left( \frac{|v(x,t)|^2}{|x|} -\frac{(v(x,t)\cdot
x)^2}{|x|^3} \right) \s'\left(\frac{|x|}{R}\right)W(|x|)\eta_\tau(t)
\,dxdt\n \\
  &&+\frac{1}{R^2}\int_{0} ^\infty\int_{\Bbb R^N}
\rho(x,t)\frac{(v(x,t)\cdot x)^2}{
|x|^2} \s^{\prime\prime} \left(\frac{|x|}{R}\right) W(|x|)\eta_\tau(t) \,dx dt\n \\
 &&+ \int_{0} ^\infty \int_{\Bbb R^N}p(x,t)\left[ W^{\prime\prime}
 (|x|) +(N-1)\frac{W' (|x|)}{|x|}\right]\sigma_R (|x|)\eta_\tau(t)
  \, dxdt\n \\
  &&+ \frac{2}{R}\int_{0} ^\infty\int_{\Bbb
R^N}p(x,t) W' (|x|)
 \s' \left(\frac{|x|}{R}\right)\eta_\tau(t)\, dxdt\n \\
&&+ \frac{N-1}{R}\int_{0} ^\infty\int_{\Bbb
R^N}p(x,t)\frac{1}{|x|}\s'
\left(\frac{|x|}{R}\right) W(|x|)\eta_\tau(t) \, dxdt\n \\
&&+ \frac{1}{R^2}\int_{0} ^\infty\int_{\Bbb R^N}p(x,t)
\s^{\prime\prime} \left(\frac{|x|}{R}\right)W(|x|)\eta_\tau(t) \,
dxdt\n \\
 &&-k\int_{0} ^\infty\int_{\Bbb R^N}
\left[\frac{|\nabla \Phi |^2}{|x|} -\frac{(x\cdot \nabla \Phi
)^2}{|x|^3} \right] W'(|x|)\si \eta_\tau(t)dxdt\n \\
 &&-k\int_{0}
^\infty\int_{\Bbb R^N} \left[\frac{|\nabla \Phi |^2}{|x|}
-\frac{(x\cdot \nabla \Phi )^2}{|x|^3}
\right]\frac{W(|x|)}{R}\sigma'
\left(\frac{|x|}{R}\right)\eta_\tau(t)\,dxdt\n \\
 &&-k\int_{0} ^\infty\int_{\Bbb R^N} \frac{(x\cdot \nabla \Phi
)^2}{|x|^2} W^{\prime\prime} (|x|)\si  \eta_\tau(t)\,dxdt\n \\
&&-k\int_{0} ^\infty\int_{\Bbb R^N} \frac{(x\cdot \nabla \Phi
)^2}{|x|^2}\left[\frac{W(|x|)}{R^2}\sigma^{\prime\prime}\left(\frac{|x|}{R}\right)
+\frac{2 W^{\prime}
(|x|)}{R}\sigma^{\prime}\left(\frac{|x|}{R}\right)\right]
\eta_\tau(t)\,dxdt
  \eqn
  \bq\label{z17}
 &&+\frac{(N-1)k}{2}\int_{0} ^\infty\int_{\Bbb R^N}
\frac{|\nabla \Phi |^2}{|x|}
W'(|x|)\si\eta_\tau(t)\,dxdt\n \\
&&+\frac{(N-1)k}{2}\int_{0} ^\infty\int_{\Bbb R^N} \frac{|\nabla
\Phi |^2}{|x|}\frac{W(|x|)}{R}\sigma' \left(\frac{|x|}{R}\right)\eta_\tau(t)\,dxdt\n \\
  &&+\frac{k}{2}\int_{0} ^\infty\int_{\Bbb R^N} |\nabla \Phi |^2 W^{\prime\prime} (|x|)\si\eta_\tau(t)\,dxdt\n \\
&&+\frac{k}{2}\int_{0} ^\infty\int_{\Bbb R^N} |\nabla \Phi |^2\left[
\frac{W(|x|)}{R^2}\sigma^{\prime\prime}\left(\frac{|x|}{R}\right)
+\frac{2 W^{\prime}
(|x|)}{R}\sigma^{\prime}\left(\frac{|x|}{R}\right)\right]\eta_\tau(t)\,dxdt\n
\\
 &&:=I_1+\cdots +I_{19}.
\eq

On the other hand, substituting $\phi(x)= \nabla\varphi_R (x)$, $\xi
(t)=\eta_\tau'(t)$
 into (\ref{1}), we find that
 \bq\label{z18}
I_3&=&\int_{0} ^\infty\int_{\Bbb R^N} \rho v(x,t)\cdot \nabla
\varphi_R (x)\eta_{\tau} '(t)\,dxdt \n \\
&=&-\int_{0} ^\infty\int_{\Bbb R^N} \rho (x,t) \varphi_R (x)
 \eta^{\prime\prime}_\tau
 (t)dxdt\n \\
 &=&-\int_{0} ^\infty\int_{\Bbb R^N} \rho (x,t)\si
 W(|x|)\eta^{\prime\prime}_\tau
 (t) dxdt\n \\
 &&\to-\int_{0} ^\infty\int_{\Bbb R^N} \rho (x,t)
W(|x|)\eta^{\prime\prime}_\tau
 (t) dxdt
 \eq
 as $R\to \infty$ by the dominated convergence theorem.
 In terms of the function $W(\cdot)$ the
condition (\ref{6})
 can be written as
\bq\label{z19} &&\int_{0} ^\infty\int_{\Bbb R^N}
(\rho(x,t)|v(x,t)|^2 +|p(x,t)|+k|\nabla \Phi (x,t) |^2)
\left[W^{\prime\prime}
 (|x|) +\right.\n \\
 &&\hspace{.5in}\left.+\frac{1}{|x|}W' (|x|) +\frac{1}{|x|^2} W(|x|)\right]
 dxdt<\infty
 \eq
 for all $T>0$.
Since
\bqn
 &&\int_{0} ^\infty\int_{\Bbb R^N} \rho(x,t)\left|\left[W^{\prime\prime} (|x|)
\frac{(v(x,t)\cdot x)^2}{|x|^2} +\right.\right.\n \\
  &&\hspace{1.5in} +\left.\left. W^{\prime} (|x|)
\left(\frac{|v(x,t)|^2}{|x|} -\frac{(v(x,t)\cdot x)^2}{|x|^3}\right)
\right]\right|\eta_\tau (t)\,dxdt\\
  &&\qquad \leq 2
 \int_{0}^{2\tau} \int_{\Bbb R^N} \rho(x,t)|v(x,t)|^2\left[ W^{\prime\prime}
(|x|)+\frac{W'(|x|)}{|x|} \right] \,dxdt <\infty, \eqn
  we can use the dominated convergence theorem to show that
  \bq\label{z20}
 && I_4 \to \int_{0} ^\infty\int_{\Bbb R^N} \rho(x,t)\left[W^{\prime\prime} (|x|)
\frac{(v(x,t)\cdot x)^2}{|x|^2} +\right. \n \\
&&\hspace{1.5in} +\left. + W^{\prime} (|x|)
\left(\frac{|v(x,t)|^2}{|x|} -\frac{(v(x,t)\cdot x)^2}{|x|^3}\right)
\right] \eta_\tau (t)\,dxdt\n \\
 \eq
  as $R\to \infty$.
Similarly,
  \bb\label{z21}
  I_8\to \int_{0} ^\infty \int_{\Bbb R^N}p(x,t)\left[ W^{\prime\prime}
 (|x|) +(N-1)\frac{W' (|x|)}{|x|}\right]\eta_\tau (t),
 \, dxdt,
\ee
 \bb
 I_{12}\to-k\int_{0} ^\infty\int_{\Bbb R^N} \left[\frac{|\nabla \Phi
|^2}{|x|} -\frac{(x\cdot \nabla \Phi )^2}{|x|^3} \right] W'(|x|)
\eta_\tau(t)dxdt,
  \ee
 \bb
 I_{14}\to -k\int_{0} ^\infty\int_{\Bbb R^N} \frac{(x\cdot \nabla \Phi
)^2}{|x|^2} W^{\prime\prime} (|x|) \eta_\tau(t)\,dxdt,
 \ee
 \bb
 I_{16}\to \frac{(N-1)k}{2}\int_{0} ^\infty\int_{\Bbb R^N}
\frac{|\nabla \Phi |^2}{|x|} W'(|x|)\eta_\tau(t)\,dxdt,
 \ee
  and
 \bb
 I_{18}\to \frac{k}{2}\int_{0} ^\infty\int_{\Bbb R^N}
 |\nabla \Phi |^2 W^{\prime\prime} (|x|)\eta_\tau(t)\,dxdt
\ee
 as $R\to \infty$.
 For $I_5$ we estimate
 \bq\label{z22}
  |I_5 |&\leq& \int_{0}^{2\tau}\int_{R<|x|<2R} \rho(x,t)|v(x,t)|^2\left|\s'
\left(\frac{|x|}{R}\right)\right|
  \frac{W'(|x|)}{|x|} \frac{|x|}{R}dxdt\n \\
  &\leq &2 \sup_{1<s<2} |\s'(s)|
\int_{0} ^{2\tau}\int_{R<|x|<2R}\rho(x) |v(x,t)|^2
  \frac{W'(|x|)}{|x|}\, dxdt\to 0 \n\\
  \eq
 as $R\to \infty$ by the dominated convergence theorem.
Similarly
  \bq\label{z23}
   |I_6|&\leq &2 \int_{0}^{2\tau}\int_{R<|x|<2R}  \frac{|x|}{R} \rho(x) |v(x,t)|^2
\left|\s'\left(\frac{|x|}{R}\right)\right|\frac{W(|x|)}{|x|^2} \,dx \n \\
&\leq &4\sup_{1<s<2} |\s'(s)|
 \int_{0}^{2\tau}\int_{R<|x|<2R}\rho(x) |v(x,t)|^2
  \frac{W'(|x|)}{|x|}\,dxdt
\to 0,\n \\
  \eq
  \bq\label{z24}
  |I_7|&\leq&\int_{0} ^{2\tau}\int_{R<|x|<2R}\frac{|x|^2}{R^2} \rho(x,t)|v(x,t)|^2\left|\s^{\prime\prime}
  \left(\frac{|x|}{R}\right)\right|\frac{W(|x|)}{|x|^2} \,dxdt\n \\
  &\leq&4\sup_{1<s<2} |\s^{\prime\prime}(s)|
 \int_{0} ^{2\tau}\int_{R<|x|<2R}\rho(x,t)|v(x,t)|^2
  \frac{W(|x|)}{|x|^2}\,dx dt\to 0,\n \\
  \eq
  and
  \bb\label{z24a}
 |I_2|\leq  2 \sup_{1<s<2}|\s ' (x)|\int_{R\leq |x|\leq 2R}
  \rho_0 (x)|v_0(x)| \frac{|W(|x|)|}{|x|} dx \to 0
  \ee
   as $R\to \infty$.
   The estimates for $I_9,I_{10}$ and $I_{11}$ are
   similar to the above, and we find
   \bq\label{z25}
   |I_9|&\leq &2 \int_{0} ^{2\tau} \int_{R<|x|<2R}|p(x,t)| \frac{|x|}{R}\frac{W'
   (|x|)}{|x|}
 \left|\s' \left(\frac{|x|}{R}\right)\right|  \, dxdt\n \\
 &\leq& 4\sup_{1<s<2} |\s' (s)|
\int_{0} ^{2\tau}\int_{R<|x|<2R}|p(x,t)|\frac{W'
   (|x|)}{|x|}\,dxdt \to 0,\n \\
  \eq
\bq\label{z26}
   |I_{10}|&\leq & (N-1)\int_{0} ^{2\tau}\int_{R<|x|<2R}|p(x,t)|\frac{|x|}{R}\left|\s'
\left(\frac{|x|}{R}\right)\right| \frac{W(|x|)}{|x|^2} \, dxdt\n \\
 &\leq& 2\sup_{1<s<2} |\s' (s)|
\int_{0} ^{2\tau}\int_{R<|x|<2R}|p(x,t)|\frac{W(|x|)}{|x|^2}dxdt\to 0,\n \\
  \eq
  and
  \bq\label{z27}
 |I_{11}|&\leq&\int_{0} ^{2\tau}\int_{\Bbb R^N}|p(x,t)|\frac{|x|^2}{R^2}
\left|\s^{\prime\prime}\left(\frac{|x|}{R}\right)\right|
\frac{W(|x|)}{|x|^2} \, dxdt \n\\
  &\leq& 4\sup_{1<s<2}
|\s^{\prime\prime} (s)|
 \int_{0} ^{2\tau}\int_{R<|x|<2R}|p(x,t)|\frac{W
   (|x|)}{|x|^2}\,dxdt \to 0\n \\
  \eq
  as $R\to \infty$ respectively.
By similar estimates we can show easily that
 \bb
  |I_{13}|+|I_{15}|+|I_{17}|+|I_{19}| \to 0.
  \ee
  as $R\to \infty$.

   Thus, passing  $R\to \infty$ in (\ref{z17}),
  we obtain
 \bq
\label{z29}
 \lefteqn{\int_{\Bbb R^N} \rho_0 (x)v_0(x)\cdot
\frac{x}{|x|}
W'(|x|)\,dx}\n \\
 &&\quad+\int_{0} ^\infty\int_{\Bbb R^N} \rho(x,t)\left[W^{\prime\prime} (|x|)
\frac{(v\cdot x)^2}{|x|^2} + W^{\prime} (|x|)
\left(\frac{|v|^2}{|x|} -\frac{(v\cdot x)^2}{|x|^3}\right)
\right]\eta_\tau (t)
 \,dxdt\n \\
  &&\quad+\int_{0} ^\infty\int_{\Bbb R^N}p(x,t)\left[
W^{\prime\prime}
 (|x|) +(N-1)\frac{W' (|x|)}{|x|}\right]\eta_\tau (t) \, dxdt\n \\
 &&\quad+\frac{\left(
 N-3\right)k}{2}\int_{0} ^\infty\int_{\Bbb R^N}\frac{|\nabla \Phi
 |^2}{|x|}W'(|x|)\eta_\tau (t)
 \,dxdt\n \\
 &&\quad+k\int_{0} ^\infty\int_{\Bbb R^N}\left[\frac{W'(|x|)}{|x|}-W^{\prime\prime}
 (|x|)\right]\frac{(x\cdot \nabla \Phi )^2}{|x|^2}\eta_\tau (t)
 \,dxdt\n \\
 &&\quad+\frac{k}{2}\int_{0} ^\infty\int_{\Bbb R^N} |\nabla \Phi |^2  W^{\prime\prime}
 (|x|)\eta_\tau (t)
 \,dxdt \n \\
 &&=\int_{0} ^\infty \int_{\Bbb R^N} \rho (x,t)
 W(|x|)\eta^{\prime\prime}_\tau .
 (t) dxdt
 \eq
 The hypothesis (\ref{5}) implies that
  \bq\label{z30}
  \lefteqn{\left|\int_{0} ^\infty \int_{\Bbb R^N} \rho (x,t)
 W(|x|)\eta^{\prime\prime}_\tau
 (t) dxdt\right|\leq \frac{1}{\tau^2} \int_{\tau}^{2\tau}
 \int_{\Bbb R^N} \rho (x,t)
 W(|x|)\left|\eta^{\prime\prime} \left(\frac{t}{\tau}\right)\right|
  dxdt }\n \\
  &\leq&\frac{1+4\tau^2}{\tau^2} \sup_{1<t<2} |\eta^{\prime\prime}(t)|
  \int_{\tau}^{2\tau}\int_{\Bbb R^N} \frac{\rho (x,t)}{1+t^2}
 \left[\int_0 ^{|x|} \int_0 ^r w(s)dsdr\right]dxdt\n\\
 &\leq& CK_1
\eq
 as $\tau\to \infty$. Next, we observe that, by our definition on
 $W(|x|)$  and the hypothesis on $w(r)$, we have
$$
W^{\prime\prime} (|x|) \frac{(v\cdot x)^2}{|x|^2} + W^{\prime} (|x|)
\left(\frac{|v|^2}{|x|} -\frac{(v\cdot x)^2}{|x|^3}\right) \geq 0,
$$
and
 $$W^{\prime\prime}
 (|x|) +(N-1)\frac{W' (|x|)}{|x|}> 0.$$
 Moreover, by the non-increasing assumption on the function $w(r)$ on $[0, \infty)$, we
 have
 \bb
 \frac{W'(|x|)}{|x|}-W^{\prime\prime}
 (|x|)=\frac{1}{|x|}\int_0 ^{|x|} w(r)dr -w(|x|)\geq 0
 \ee
 for almost every $x\in \Bbb R^N$.
Thus, we can apply the monotone convergence theorem to obtain
 \bq\label{z30a}
&&\int_{0} ^\infty\int_{\Bbb R^N} \rho(x,t)\left[W^{\prime\prime}
(|x|) \frac{(v\cdot x)^2}{|x|^2} + W^{\prime} (|x|)
\left(\frac{|v|^2}{|x|} -\frac{(v\cdot x)^2}{|x|^3}\right)
\right]\eta_\tau (t)
 \,dxdt\n \\
 &&\quad\to \int_{0} ^\infty\int_{\Bbb R^N} \rho(x,t)\left[W^{\prime\prime}
(|x|) \frac{(v\cdot x)^2}{|x|^2} + W^{\prime} (|x|)
\left(\frac{|v|^2}{|x|} -\frac{(v\cdot x)^2}{|x|^3}\right) \right]
 \,dxdt,\n \\
 \eq
  \bq\label{z30b}
 &&\int_{0} ^\infty\int_{\Bbb R^N}p(x,t)\left[
W^{\prime\prime}
 (|x|) +(N-1)\frac{W' (|x|)}{|x|}\right]\eta_\tau (t) \, dxdt\n \\
 &&\quad \to \int_{0} ^\infty\int_{\Bbb R^N} \rho(x,t)\left[W^{\prime\prime}
(|x|) \frac{(v\cdot x)^2}{|x|^2} + W^{\prime} (|x|)
\left(\frac{|v|^2}{|x|} -\frac{(v\cdot x)^2}{|x|^3}\right) \right]
 \,dxdt,\n \\
 \eq
\bq
 \int_{0} ^\infty\int_{\Bbb R^N}\frac{|\nabla \Phi
 |^2}{|x|}W'(|x|)\eta_\tau (t)
 \,dxdt \to\int_{0} ^\infty\int_{\Bbb R^N}\frac{|\nabla \Phi
 |^2}{|x|}W'(|x|)
 \,dxdt
 \eq
 \bq
 &&\int_{0} ^\infty\int_{\Bbb R^N}\left[\frac{W'(|x|)}{|x|}-W^{\prime\prime}
 (|x|)\right]\frac{(x\cdot \nabla \Phi )^2}{|x|^2}\eta_\tau (t)
 \,dxdt\n \\
 &&\qquad \qquad \to \int_{0} ^\infty\int_{\Bbb R^N}\left[\frac{W'(|x|)}{|x|}-W^{\prime\prime}
 (|x|)\right]\frac{(x\cdot \nabla \Phi )^2}{|x|^2}
 \,dxdt,
 \eq
 \bq
&&\int_{0} ^\infty\int_{\Bbb R^N} |\nabla \Phi |^2 W^{\prime\prime}
 (|x|)\eta_\tau (t)
 \,dxdt\to\int_{0} ^\infty\int_{\Bbb R^N} |\nabla \Phi |^2
W^{\prime\prime}
 (|x|)
 \,dxdt\n \\
 \eq
 as $\tau \to \infty$. Thus, passing $\tau\to \infty$ in (\ref{z29}), we
 find that
\bqn
 \lefteqn{\int_{\Bbb R^N} \rho_0 (x)v_0(x)\cdot
\frac{x}{|x|}
W'(|x|)\,dx}\n \\
 &&+\int_{0} ^\infty\int_{\Bbb R^N} \rho(x,t)\left[W^{\prime\prime} (|x|)
\frac{(v\cdot x)^2}{|x|^2} + W^{\prime} (|x|)
\left(\frac{|v|^2}{|x|} -\frac{(v\cdot x)^2}{|x|^3}\right) \right]
 \,dxdt\n \\
&&+\int_{0} ^\infty\int_{\Bbb R^N}p(x,t)\left[ W^{\prime\prime}
 (|x|) +(N-1)\frac{W' (|x|)}{|x|}\right] \, dxdt
 \eqn
\bq\label{z30c} &&+\frac{\left(
 N-3\right)k}{2}\int_{0} ^\infty\int_{\Bbb R^N}\frac{|\nabla \Phi
 |^2}{|x|}W'(|x|)
 \,dxdt\n \\
 &&+k\int_{0} ^\infty\int_{\Bbb R^N}\left[\frac{W'(|x|)}{|x|}-W^{\prime\prime}
 (|x|)\right]\frac{(x\cdot \nabla \Phi )^2}{|x|^2}
 \,dxdt\n \\
 &&+\frac{k}{2}\int_{0} ^\infty\int_{\Bbb R^N} |\nabla \Phi |^2  W^{\prime\prime}
 (|x|)
 \,dxdt \leq CK_1 ,
 \eq
 which proves (\ref{6a}) for $N\geq 1$.  In the case $N=2$, we
 choose $w(r)\equiv 1$ on $[0, \infty)$. Then, the inequality
 (\ref{z30c}) reduces to
\bb\label{z30ca} \int_{\Bbb R^2} \rho_0 (x)v_0(x)\cdot x\,dx
+\int_{0} ^\infty\int_{\Bbb R^2} [\rho(x,t)|v(x,t)|^2 +2p(x,t)]
 \,dxdt \leq CK_1.
 \ee
 $\square$\\
 \ \\
 In order to establish Theorem 1.2 we use the following lemma, which is proved in \cite{guo}.
\begin{lemma}
Suppose $(\rho, v)$ is a global weak solution of $(NS)$ with the
setting given by Theorem 1.2.  We suppose that the energy inequality
(\ref{energy}) holds. Then,
 \bb\label{guo}
 \int_0 ^\infty \int_{\Bbb R^N}
 \frac{ \rho(x,t)(1+|x|^2)^{\frac{N+2}{4\gamma}}}{t^2} dxdt \leq C E(0).
 \ee
\end{lemma}
Since $\frac{N+2}{4\gamma}\geq 1$ in our setting of Theorem 1.2, one
immediate consequence of (\ref{guo}) is the following fact
  \bb\label{key}
\lim_{\tau\to \infty} \int_{\tau}^{2\tau}\int_{\Bbb R^N} \frac{\rho
(x,t)}{1+t^2} |x|^2dxdt=0.
 \ee
 Indeed, using (\ref{guo}), we deduce
 $$
\lim_{\tau\to \infty} \int_{\tau}^{2\tau}\int_{\Bbb R^N} \frac{\rho
(x,t)}{1+t^2} |x|^2dxdt\leq\lim_{\tau\to \infty}
\int_{\tau}^{2\tau}\int_{\Bbb R^N}
 \frac{ \rho(x,t)(1+|x|^2)^{\frac{N+2}{4\gamma}}}{t^2} dxdt=0,
 $$
 where the last step follows from the dominated convergence
 theorem.\\
 \ \\
 \noindent{\bf Proof of Theorem 1.3 }
Suppose there exists a global weak solution $(\rho, v, \Phi)$
satisfying (\ref{1})-(\ref{4})(with $\mu \neq 0$). Here, we choose
the vector test function as
 \bb\label{z33}
\phi=\nabla \varphi_R, \quad\varphi_R (x)=\frac{|x|^2}{2} \s
\left(\frac{|x|}{R}\right)=\frac{|x|^2}{2}\s_R (|x|),
 \ee
 where $\s$ is the cut-off function defined in (\ref{18}).
 Similarly to the proof of  Theorem 1.2 we also use the same temporal cut-off
 function $\eta_\tau(t)$ defined in (\ref{z15})-(\ref{z16}). Substituting
$\phi(x)=\nabla \varphi_R (x), \xi(t)=\eta_\tau(t)$ into (\ref{2}),
we have
  \bqn
 \lefteqn{0=\int_{\Bbb R^N} \rho_0 (x)v_0(x)\cdot x\si dx+\frac{1}{2 R}\int_{\Bbb R^N}
  \rho_0 (x)v_0(x)\cdot x |x|\s'\left(\frac{|x|}{R}\right)
 dx}\hspace{.0in}\n \\
 &&+\int_0 ^\infty \int_{\Bbb R^N} \rho (x,t)v(x,t) \cdot
 \nabla \varphi_R (x)
 \eta_\tau ^{\prime}(t) dxdt \n \\
&&+\int_0 ^\infty\int_{\Bbb R^N} \rho (x,t) |v(x,t)|^2 \si \eta_\tau(t) \,dx dt\n\\
&&+\frac{1}{2R}\int_0 ^\infty\int_{\Bbb R^N} \rho(x,t)  \s'
\left(\frac{|x|}{R}\right) \frac{(v(x,t)\cdot x)^2}{|x|}\eta_\tau(t)
\,dxdt\n \\
 &&+ \frac{1}{2R}\int_0 ^\infty\int_{\Bbb R^N}
\rho(x,t)|v(x,t)|^2
|x| \s'\left(\frac{|x|}{R}\right)\eta_\tau(t) \,dxdt \n \\
  &&+\frac{1}{2R^2}\int_0 ^\infty\int_{\Bbb R^N} \rho(x,t)(v(x,t)\cdot x)^2}{
\s^{\prime\prime} \left(\frac{|x|}{R}\right) \eta_\tau(t) \,dx dt\n \\
 &&+ N\int_0 ^\infty \int_{\Bbb R^N}p(x,t)\sigma_R (|x|)\eta_\tau(t)
  \, dxdt\n \\
&&+ \frac{2}{R}\int_0 ^\infty\int_{\Bbb R^N}p(x,t) |x|
 \s' \left(\frac{|x|}{R}\right)\eta_\tau(t)\, dxdt\n \\
&&+ \frac{N-1}{2R}\int_0 ^\infty\int_{\Bbb R^N}p(x,t)|x|\s'
\left(\frac{|x|}{R}\right) \eta_\tau(t) \, dxdt\n \\
&&+ \frac{1}{2R^2}\int_0 ^\infty\int_{\Bbb R^N}p(x,t)|x|^2
\s^{\prime\prime} \left(\frac{|x|}{R}\right)\eta_\tau(t) \, dxdt\n
\\
&&-k\int_{0} ^\infty\int_{\Bbb R^N} \left[|\nabla \Phi |^2
-\frac{(x\cdot \nabla \Phi )^2}{|x|^2} \right] \si
\eta_\tau(t)dxdt\n \\
&&-k\int_{0} ^\infty\int_{\Bbb R^N} \left[|\nabla \Phi |^2
-\frac{(x\cdot \nabla \Phi )^2}{|x|^2} \right]\frac{|x|}{2R}\sigma'
\left(\frac{|x|}{R}\right)\eta_\tau(t)\,dxdt\n \\
 &&-k\int_{0} ^\infty\int_{\Bbb R^N} \frac{(x\cdot \nabla \Phi
)^2}{|x|^2} \si  \eta_\tau(t)\,dxdt\n \\
&&-k\int_{0} ^\infty\int_{\Bbb R^N}(x\cdot \nabla \Phi
)^2\left[\frac{1}{2R^2}\sigma^{\prime\prime}\left(\frac{|x|}{R}\right)
+\frac{2}{|x|R}\sigma^{\prime}\left(\frac{|x|}{R}\right)\right]
\eta_\tau(t)\,dxdt\n \\
 &&+\frac{(N-1)k}{2}\int_{0} ^\infty\int_{\Bbb R^N}
|\nabla \Phi |^2 \si\eta_\tau(t)\,dxdt
  \eqn
  \bq\label{z36}
&&+\frac{(N-1)k}{4}\int_{0} ^\infty\int_{\Bbb R^N}|\nabla \Phi
|^2\frac{|x|}{R}\sigma'
\left(\frac{|x|}{R}\right)\eta_\tau(t)\,dxdt\n \\
  &&+\frac{k}{2}\int_{0} ^\infty\int_{\Bbb R^N} |\nabla \Phi |^2 \si\eta_\tau(t)\,dxdt\n \\
&&+\frac{k}{2}\int_{0} ^\infty\int_{\Bbb R^N} |\nabla \Phi |^2\left[
\frac{|x|^2}{2R^2}\sigma^{\prime\prime}\left(\frac{|x|}{R}\right)
+\frac{2
|x|}{R}\sigma^{\prime}\left(\frac{|x|}{R}\right)\right]\eta_\tau(t)\,dxdt\n
\\
&&+(2\mu +\lambda)\int_0 ^\infty\int_{\Bbb R^N} v\cdot\nabla \Delta
(|x|^2
 \s \left(\frac{|x|}{R}\right)\eta_\tau(t)\, dxdt,\n \\
  &&:=I_1+\cdots +I_{20}.
 \eq
On the other hand, substituting $\phi(x)= \nabla\varphi_R (x)$, $\xi
(t)=\eta_\tau'(t)$
 into (\ref{1}), then similarly as before, we find that (note that $\xi (0)=\eta_\tau'(0)=0$)
 \bq\label{z37}
 I_3&=&\int_0 ^\infty\int_{\Bbb R^N} \rho v(x,t)\cdot
\nabla
\varphi_R (x)\eta_{\tau} '(t)\,dxdt\n \\
&=&- \int_0 ^\infty \int_{\Bbb R^N} \rho (x,t) \varphi_R (x)
 \eta^{\prime\prime}_\tau
 (t)dxdt\n \\
 &=&-\int_0 ^\infty \int_{\Bbb R^N} \rho (x,t)|x|^2\si \eta^{\prime\prime}_\tau
 (t) dxdt\n \\
 &&\to-\int_0 ^\infty \int_{\Bbb R^N} \rho (x,t)
|x|^2\eta^{\prime\prime}_\tau
 (t) dxdt
 \eq
 as $R\to \infty$ by the dominated convergence theorem.
  We  also have
  \bb\label{z38}
  I_4 \to \int_0 ^\infty\int_{\Bbb R^N} \rho(x,t)|v(x,t)|^2  \eta_\tau (t)\,dxdt
 \ee
  as $R\to \infty$.
Similarly,
 \bb\label{z38a}
 I_1\to \int_{\Bbb R^N} \rho_0 (x)v_0(x)\cdot x dx,
 \ee
  \bb\label{z39}
  I_8\to N\int_0 ^\infty \int_{\Bbb R^N}p(x,t)\eta_\tau (t)
 \, dxdt,
\ee
 \bb
I_{12}\to-k\int_{0} ^\infty\int_{\Bbb R^N} \left[|\nabla \Phi |^2
-\frac{(x\cdot \nabla \Phi )^2}{|x|^2} \right]  \eta_\tau(t)dxdt,
 \ee
 \bb
I_{14}\to-k\int_{0} ^\infty\int_{\Bbb R^N} \frac{(x\cdot \nabla \Phi
)^2}{|x|^2}  \eta_\tau(t)\,dxdt,
 \ee
 \bb
I_{16}\to\frac{(N-1)k}{2}\int_{0} ^\infty\int_{\Bbb R^N} |\nabla
\Phi |^2 \eta_\tau(t)\,dxdt,
 \ee
 and
 \bb
 I_{18}\to\frac{k}{2}\int_{0} ^\infty\int_{\Bbb R^N} |\nabla \Phi |^2\eta_\tau(t)\,dxdt
 \ee
  as $R\to \infty$.
 For $I_5, I_6$ we estimate
 \bq\label{z40}
  |I_5 |+|I_6|&\leq& \int_0 ^{2\tau}\int_{R<|x|<2R} \rho(x,t)|v(x,t)|^2\left|\s'
\left(\frac{|x|}{R}\right)\right|
  \frac{|x|}{R}dxdt\n \\
  &\leq &2 \sup_{1<s<2} |\s'(s)|
\int_0 ^{2\tau}\int_{R<|x|<2R}\rho(x) |v(x,t)|^2
 \, dxdt
\to 0\n \\
  \eq
 as $R\to \infty$ by the dominated convergence theorem.
Similarly
 \bb\label{z40a}
  |I_2|\leq \int_{R<|x|<2R} \rho_0 (x) |x| dx \to 0,
\ee
 and
 \bq\label{z41}
  |I_7|&\leq&\frac12\int_0 ^{2\tau}\int_{R<|x|<2R}\frac{|x|^2}{R^2}
   \rho(x,t)|v(x,t)|^2\left|\s^{\prime\prime}
  \left(\frac{|x|}{R}\right)\right| \,dxdt\n \\
  &\leq&2\sup_{1<s<2} |\s^{\prime\prime}(s)|
 \int_0 ^{2\tau}\int_{R<|x|<2R}\rho(x,t)|v(x,t)|^2
 \,dx \to 0\n \\
  \eq
   as $R\to \infty$. The estimates for $I_9,I_{10}$ and $I_{11}$ are
   similar to the above, and we find
   \bq\label{z42}
   |I_9|&\leq &2 \int_0 ^{2\tau} \int_{R<|x|<2R}|p(x,t)| \frac{|x|}{R}
 \left|\s' \left(\frac{|x|}{R}\right)\right|  \, dxdt\n \\
 &\leq& 4\sup_{1<s<2} |\s' (s)|
\int_0 ^{2\tau}\int_{R<|x|<2R}|p(x,t)|\,dxdt \to 0,\n \\
  \eq
\bq\label{z43}
   |I_{10}|&\leq & \frac{N-1}{2R}\int_0 ^{2\tau}\int_{R<|x|<2R}|p(x,t)||x|\left|\s'
\left(\frac{|x|}{R}\right)\right|  \, dxdt\n \\
 &\leq& (N-1)\sup_{1<s<2} |\s' (s)|
\int_0 ^{2\tau}\int_{R<|x|<2R}|p(x,t)|dxdt\to 0,\n \\
  \eq
  and
  \bq\label{z44}
 |I_{11}|&\leq&\frac{1}{2R^2} \int_0 ^{2\tau}\int_{\Bbb
 R^N}|p(x,t)||x|^2
\left|\s^{\prime\prime}\left(\frac{|x|}{R}\right)\right|
 \, dxdt \n\\
  &\leq& 2\sup_{1<s<2}
|\s^{\prime\prime} (s)|
 \int_0 ^{2\tau}\int_{R<|x|<2R}|p(x,t)|\,dxdt \to 0\n \\
  \eq
  as $R\to \infty$ respectively. Similarly,
  \bb
  |I_{13}|+|I_{15}|+|I_{17}|+|I_{19}|\to 0
 \ee
 as $R\to \infty$.  Now we show the vanishing of the
  viscosity term as $R\to \infty$.
  Similarly to stationary case we estimate
  \bq\label{z45}
  |I_{20}|&=&(2\mu +\lambda ) \left|\int_0 ^\infty\int_{\Bbb R^N} v\cdot\nabla \Delta (|x|^2
 \s \left(\frac{|x|}{R}\right)\eta_\tau(t)\,dxdt\right|\n \\
&\leq& (2\mu +\lambda) \left|\int_0 ^\infty\int_{\Bbb R^N}
(N+5)\left[\frac{(v\cdot x)}{R|x|}\s ' \left(\frac{|x|}{R}\right)+
\frac{(v\cdot x)}{R^2}\s^{\prime\prime}
\left(\frac{|x|}{R}\right)\right]\eta_\tau(t)\, dxdt\right|\n \\
&&\qquad+(2\mu +\lambda)\left|\int_0 ^\infty\int_{\Bbb R^N}
\frac{|x|(v\cdot x)}{R^3}\s^{\prime\prime\prime}
\left(\frac{|x|}{R}\right)\eta_\tau(t)\,dx dt\right| \n \\
 &\leq&
\frac{C}{R} \int_0 ^{2\tau}\int_{R\leq |x|\leq 2R} |v(x,t)|
\, dx dt\n \\
&\leq& C \int_0 ^{2\tau}\left(\int_{R\leq  |x|\leq 2R}
|v(x)|^{\frac{N}{N-1}} dx\right)^{\frac{N-1}{N}} dt\to 0
 \eq
   as $R\to \infty$.
   Thus passing  $R\to \infty$ in (\ref{z36}),
  we obtain
\bq\label{z46}
 &&\int_0 ^\infty \int_{\Bbb R^N} \rho (x,t)
|x|^2\eta^{\prime\prime}_\tau
 (t) dxdt=\int_0 ^\infty\int_{\Bbb R^N} \rho(x,t) |v(x,t)|^2\eta_\tau (t)
 \,dxdt \n \\
&&\qquad+N\int_0 ^\infty\int_{\Bbb R^N}p(x,t)\eta_\tau
(t)dxdt+\frac{N-2}{2}\int_{0} ^\infty\int_{\Bbb R^N} |\nabla \Phi
|^2 \eta_\tau(t)\,dxdt\n \\
&&\qquad+ \int_{\Bbb R^N} \rho_0 (x)v_0(x)\cdot x
 dx
 \eq
 for any $\tau >0$.
Note that
  \bq\label{z48}
   \lefteqn{\left|\int_0 ^\infty \int_{\Bbb R^N} \rho (x,t)
|x|^2\eta^{\prime\prime}_\tau
 (t) dxdt\right|\leq \frac{1}{\tau^2} \int_{\tau}^{2\tau}
 \int_{\Bbb R^N} \rho (x,t)
|x|^2\left|\eta^{\prime\prime} \left(\frac{t}{\tau}\right)\right|
  dxdt }\n \\
  &&\leq\frac{1+4\tau^2}{\tau^2} \sup_{1<s<2} |\eta^{\prime\prime}(s)|
  \int_{\tau}^{2\tau}\int_{\Bbb R^N} \frac{\rho (x,t)}{1+t^2}
|x|^2dxdt\to 0 ,
 \eq
  as $\tau\to \infty$ by (\ref{key}). On the other hand, by the monotone convergence theorem we deduce
 \bq
&&\int_0 ^\infty\int_{\Bbb R^N} \rho(x,t) |v(x,t)|^2\eta_\tau (t)
\,dxdt \to \int_0 ^\infty\int_{\Bbb R^N} \rho(x,t) |v(x,t)|^2
\,dxdt,\n \\
&&\int_0 ^\infty\int_{\Bbb R^N}p(x,t)\eta_\tau (t)\,dxdt\to \int_0
^\infty\int_{\Bbb R^N}p(x,t)\,dxdt\n \\
&&\int_{0} ^\infty\int_{\Bbb R^N} |\nabla \Phi |^2
\eta_\tau(t)\,dxdt\to \int_{0} ^\infty\int_{\Bbb R^N} |\nabla \Phi
|^2\,dxdt \eq
 as $\tau\to \infty$. Thus, passing $\tau \to \infty$ in
(\ref{z46})  we
 have
\bq\label{z48a}
 \lefteqn{0=\int_0 ^\infty\int_{\Bbb R^N} \rho(x,t) |v(x,t)|^2
 \,dxdt +N\int_0 ^\infty\int_{\Bbb R^N}p(x,t)\,dxdt}\hspace{.0in}\n \\
&&+ \frac{N-2}{2}\int_{0} ^\infty\int_{\Bbb R^N} |\nabla \Phi
|^2\,dxdt +\int_{\Bbb R^N} \rho_0 (x)v_0(x)\cdot x\,
 dx.
 \eq
$\square$\\

$$\mbox{ \bf Acknowledgements} $$
 This
work was supported partially by  KRF Grant(MOEHRD, Basic Research
Promotion Fund).

\end{document}